\documentclass[11pt,letterpaper,reqno]{amsart}
\usepackage{amsmath,amsthm,amsfonts,amssymb,mathtools,bbm}
\usepackage[lite]{amsrefs}

\newtheorem{theorem}{Theorem}

\newtheorem{corollary}[theorem]{Corollary}
\newtheorem{conjecture}[theorem]{Conjecture}

\numberwithin{equation}{section}

\begin{document}

\title{On a trilinear singular integral form with determinantal kernel}

\author[P. Gressman]{Philip Gressman}
\address{Philip Gressman, University of Pennsylvania, Department of Mathematics, Dav
id Rittenhouse Lab, 209
South 33rd Street, Philadelphia, PA 19104-6395, USA}
\email{gressman@math.upenn.edu}

\author[D. He]{Danqing He}
\address{Danqing He, University of Missouri, Department of Mathematics, Columbia, Missouri}
\email{dhd27@mail.missouri.edu}

\author[V. Kova\v{c}]{Vjekoslav Kova\v{c}}
\address{Vjekoslav Kova\v{c}, University of Zagreb, Faculty of Science, Department of Mathematics, Bijeni\v{c}ka cesta 30, 10000 Zagreb, Croatia}
\email{vjekovac@math.hr}

\author[B. Street]{Brian Street}
\address{Brian Street, University of Wisconsin-Madison, Department of Mathematics, 480 Lincoln Dr., Madison, WI, 53706}
\email{street@math.wisc.edu}

\author[C. Thiele]{Christoph Thiele}
\address{Christoph Thiele, Universit\"at Bonn, Endenicher Allee 60, 53115 Bonn, Germany}
\email{thiele@math.uni-bonn.de}

\author[P.-L. Yung]{Po-Lam Yung}
\address{Po-Lam Yung, The Chinese University of Hong Kong, Ma Liu Shui, Hong Kong}
\email{plyung@math.cuhk.edu.hk}



\begin{abstract}
We study a trilinear singular integral form acting on two-dimensional functions and possessing invariances under arbitrary matrix dilations and linear modulations. One part of the motivation for introducing it lies in its large symmetry groups acting on the Fourier side. Another part of the motivation is that this form stands between the bilinear Hilbert transforms and the first Calder\'{o}n commutator, in the sense that it can be reduced to a superposition of the former, while it also successfully encodes the latter. As the main result we determine the exact range of exponents in which the $\textup{L}^p$ estimates hold for the considered form.
\end{abstract}

\date{\today}

\maketitle

\section{Introduction}

This paper studies a highly symmetric trilinear form, which is a singular integral variant of the determinant functionals considered in \cite{Gre}. Related singular integral variants have been studied in \cite{Val}. We refer to \cite{Gre} and \cite{Val} for a discussion of prior work in this context. For three Schwartz functions $f,g,h$ defined on $\mathbb{R}^2$ we define the trilinear form
\begin{align}
\Lambda(f,g,h) & := \textup{p.v.} \int_{(\mathbb{R}^2)^3} f(x)g(y)h(z)\delta(x+y+z)
\det\begin{pmatrix}1 & 1 & 1\\ x & y & z \end{pmatrix}^{-1} dxdydz \label{lambdaspace} \\
& \,= \frac{1}{3} \,\textup{p.v.} \int_{(\mathbb{R}^2)^3} \widehat{f}(\xi)\widehat{g}(\eta)\widehat{h}(\zeta)
\det\begin{pmatrix}1 & 1 & 1\\ \xi & \eta & \zeta \end{pmatrix}^{-1} d\xi d\eta d\zeta. \label{lambdafrequency}
\end{align}
Note that $x,y,z$ are elements of $\mathbb{R}^2$, written as $2\times 1$ columns, and thus the determinant in \eqref{lambdaspace} is taken of a $3\times 3$ matrix. The principal value is defined as the limit as $\epsilon\to 0$ of the integral over the set of triples $(x,y,z)$ for which the absolute value of the determinant in \eqref{lambdaspace} exceeds $\epsilon$. More details on the definition, in particular the Dirac delta notation, as well as justification of the equality between \eqref{lambdaspace} and \eqref{lambdafrequency}, will be presented in Section~\ref{convergence}.

The main result of the present paper is:
\begin{theorem}\label{triltheorem}
If
\begin{equation}\label{eqexpreg}
2<p,q,r<\infty, \quad \frac{1}{p} + \frac{1}{q} + \frac{1}{r} = 1,
\end{equation}
then
\begin{equation}\label{eqformbound1}
|\Lambda(f,g,h)| \lesssim_{p,q,r} \|f\|_{\textup{L}^{p}(\mathbb{R}^2)} \|g\|_{\textup{L}^{q}(\mathbb{R}^2)} \|h\|_{\textup{L}^{r}(\mathbb{R}^2)}.
\end{equation}
The range of exponents is sharp in the sense that for $1\leq p,q,r\leq \infty$ not satisfying \eqref{eqexpreg} the a priori estimate \eqref{eqformbound1} fails.
\end{theorem}
Here and in what follows we write $A\lesssim_P B$ for two nonnegative quantities $A$ and $B$ whenever $A\leq C_P B$ holds with some finite constant $C_P$ depending on a set of parameters $P$.

Combining \eqref{eqformbound1} with the Hausdorff-Young inequality applied to the three functions $\widehat{f},\widehat{g},\widehat{h}$ we obtain the following considerably simpler corollary, which appears in Valdimarsson's work \cite{Val}:
\begin{corollary}[\cite{Val}]\label{valdicorollary}
Let $p,q,r$ be as in \eqref{eqexpreg} and let $p',q',r'$ denote the dual exponents. Then
\begin{equation}
|\Lambda(f,g,h)| \lesssim_{p,q,r} \|\widehat{f}\|_{\textup{L}^{p'}(\mathbb{R}^2)} \|\widehat{g}\|_{\textup{L}^{q'}(\mathbb{R}^2)} \|\widehat{h}\|_{\textup{L}^{r'}(\mathbb{R}^2)}.
\end{equation}
\end{corollary}
It is also shown in \cite{Val} that this range of exponents is optimal, and certain endpoint bounds at its boundary are discussed.

We present the proof of Theorem~\ref{triltheorem} in Section~\ref{mainproof}. The key elements are a non-linear change of variables followed by a fiber-wise application of uniform estimates for the bilinear Hilbert transform. The latter, in its dual form, is a trilinear form on single-variable functions $\tilde{f},\tilde{g},\tilde{h}$ defined by
\begin{equation}\label{bht}
\mathcal{B}(\tilde{f},\tilde{g},\tilde{h}) := \textup{p.v.} \int_{\mathbb{R}^3} \tilde{f}(x)\tilde{g}(y)\tilde{h}(z) \delta(ax+by+cz)
\det\begin{pmatrix}1 & 1 & 1\\ a^{-1} & b^{-1} & c^{-1} \\x & y & z  \end{pmatrix}^{-1} dxdydz
\end{equation}
for any fixed unit vector $(a,b,c)$ perpendicular to $(1,1,1)$ whose components are all non-zero; see Section~\ref{mainproof} for the more common formulation of this definition. Bounds analogous to Theorem~\ref{triltheorem} for this operator for a fixed triple $(a,b,c)$ were established in \cite{LT1}. Limiting forms as either of the components of $(a,b,c)$ tends to zero reduce by renormalization to the linear Hilbert transform in combination with a pointwise product, so these limit forms satisfy analogous estimates with much simpler proofs. Uniform bounds on \eqref{bht} as one of the components approaches $0$ were first considered in \cite{Thi} and settled in the range required for proving Theorem~\ref{triltheorem} by Grafakos and Li in \cite{GL}. We note that the exact range of exponents for which one has bounds for \eqref{bht} is larger than that in Theorem~\ref{triltheorem}; see \cite{LT2} and \cite{Li}. It is a rare feat in the context of time-frequency analysis that one understands the exact range of exponents as in Theorem~\ref{triltheorem}.

Proving bounds for some trilinear form by means of superposition of bilinear Hilbert transforms \eqref{bht} has a prominent precedent. That way one can establish estimates for the so-called first commutator of Calder\'{o}n. Indeed, the idea of writing the commutator as a superposition of bilinear Hilbert transforms was purportedly one of Calder\'{o}n's motivations for considering the bilinear Hilbert transform in the 1960s; see \cite{Cal}. The history, of course, is that he proceeded to prove bounds for his first commutator by different means, leaving the more difficult bilinear Hilbert transform open.  For Calder\'{o}n's suggested approach bounds for the bilinear Hilbert transform with dependence
$$ O(|a|^{-\epsilon}+|b|^{-\epsilon}+|c|^{-\epsilon}) $$
in the parameter $(a,b,c)$ for some $0<\epsilon<1$ are sufficient. A distinctive feature of our present application is that such dependence in $\epsilon$ for any fixed $\epsilon>0$ appears to be insufficient to derive Theorem~\ref{triltheorem}, at least with our superposition argument in the full range of exponents. In this sense our present application is a much sharper use of the uniform bounds for the bilinear Hilbert transform. In Section~\ref{remarks} we will further manifest that \eqref{lambdaspace} lies between the bilinear Hilbert transform and Calder\'{o}n's commutator by presenting an argument that Theorem~\ref{triltheorem} implies the corresponding bounds for Calder\'{o}n's commutator.

The form \eqref{lambdaspace} and its frequency domain counterpart \eqref{lambdafrequency} are invariant under actions of several classical groups; we refer to \cite{Oks},\cite{Pra} for a hint of importance of invariant trilinear forms in representation theory. Let us first observe that $\Lambda$ possesses the full $\textup{GL}_2(\mathbb{R})$ dilation symmetry. Namely, if for a nonsingular $2\times 2$ real matrix $A$ we denote
$$ (\textup{D}_A f)(x) := (\mathop{\textup{sgn}}(\det A)) |\det A|^{-1/3} f(A^{-1}x), $$
then a simple linear change of variables gives
\begin{equation}\label{eqdilatinv}
\Lambda\big(\textup{D}_A f, \textup{D}_A g, \textup{D}_A h\big) = \Lambda(f,g,h).
\end{equation}
Observe that $\Lambda$ is not invariant under translations, but it has the most obvious modulation invariance: if for $b\in\mathbb{R}^2$ we set
$$ (\textup{M}_b f)(x) := e^{2\pi i b\cdot x} f(x), $$
where ``$\cdot$'' denotes the standard inner product in $\mathbb{R}^2$, then
\begin{equation}\label{eqmodulinv}
\Lambda(\textup{M}_b f, \textup{M}_b g, \textup{M}_b h) = \Lambda(f,g,h).
\end{equation}
Yet another symmetry of $\Lambda$ is the alternating property under the permutations of the functions:
$$ \Lambda(f,g,h) = \Lambda(g,h,f) = \Lambda(h,f,g) = -\Lambda(g,f,h) = -\Lambda(h,g,f) = -\Lambda(f,h,g). $$

We can combine \eqref{eqdilatinv} and \eqref{eqmodulinv} into a single invariance of \eqref{lambdafrequency} under simultaneous action of the affine group $\textup{Aff}(\mathbb{R}^2)$ on the Fourier transforms of the three functions. The affine group in the plane consists of all transformations of the form $\xi\mapsto L\xi+\rho$ for an invertible linear map $L$ on $\mathbb{R}^2$ and a vector $\rho\in\mathbb{R}^2$. If $\tilde{f},\tilde{g},\tilde{h}$ are related to $f,g,h$ via
$$ \widehat{\tilde{f}}(\xi) = (\mathop{\textup{sgn}}(\det L)) |\det L|^{2/3} \widehat{f}(L\xi+\rho) $$
and the other two analogous relations, then we have
$$ \Lambda(\tilde{f},\tilde{g},\tilde{h}) = \Lambda(f,g,h). $$
The affine group in the plane is six dimensional and maps transitively on the set of triples of points $(\xi,\eta,\zeta)\in (\mathbb{R}^2)^3$ that are not collinear. Therefore, the functions on $(\mathbb{R}^2)^3$ invariant under this action are determined up to a scalar multiple on a Zariski open set. 

Finally, (\ref{eqdilatinv}) can be viewed as part of a bigger symmetry of $\Lambda$. Consider $\mathbb{R}^2$ as a Zariski open subset of the real projective plane $\mathbb{R}\textup{P}^2$. The general linear group $\textup{GL}_3(\mathbb{R})$ acts on $\mathbb{R}\textup{P}^2$, and this action induces an invariance of our trilinear form $\Lambda$ under $\textup{GL}_3(\mathbb{R})$. More precisely, suppose $M \in \textup{GL}_3(\mathbb{R})$. For $\xi = (\xi_1,\xi_2) \in \mathbb{R}^2$ which we identify with the equivalent class $[1,\xi_1,\xi_2] \in \mathbb{R}\textup{P}^2$, let $\tilde{\xi} = (\tilde{\xi}_0, \tilde{\xi}_1, \tilde{\xi}_2)$ be the vector in $\mathbb{R}^3$ given by $\tilde{\xi}^t = M (1,\xi_1,\xi_2)^t$, and let $(\breve{\xi}_1,\breve{\xi}_2) \in \mathbb{R}^2$ be such that $\tilde{\xi} \in [1,\breve{\xi}_1,\breve{\xi}_2]$. Define $\breve{f}$ by
$$\widehat{\breve{f}}(\xi) = (\mathop{\textup{sgn}} (\det M)) |\det M|^{-1/3} \left| \det \left( \frac{\partial \breve{\xi}}{\partial \xi} \right) \right| \tilde{\xi}_0 \widehat{f}(\breve{\xi})$$ and similarly $\breve{g}$ and $\breve{h}$. Then $$\Lambda(\breve{f},\breve{g},\breve{h}) = \Lambda(f,g,h)$$ for any $M \in \textup{GL}_3(\mathbb{R})$. In particular, specifying to $M$ of the form $\left(
\begin{array}{cc}
1 & 0 \\ 
0 & L 
\end{array}
\right)$ with $L \in \textup{GL}_2(\mathbb{R})$, we recover invariance (\ref{eqdilatinv}); specifying to the rotation matrix
$
\left(
\begin{array}{ccc}
0 & 1 & 0 \\
-1 & 0 & 0 \\
0 & 0 & 1 
\end{array}
\right),
$
we get $\Lambda(\breve{f},\breve{g},\breve{h}) = \Lambda(f,g,h)$ if 
$$
\widehat{\breve{f}}(\xi_1,\xi_2) = (\mathop{\textup{sgn}}\xi_1)\xi_1^{-2}\widehat{f}(-\xi_1^{-1},\xi_1^{-1}\xi_2)
$$
and similarly for $\breve{g}$, $\breve{h}$.

Valdimarsson \cite{Val} also proves a multilinear version of Corollary~\ref{valdicorollary} with arbitrary degree of multilinearity. Consequently, it is natural to conjecture that an analogous multilinear generalization of Theorem~\ref{triltheorem} may hold. We formulate this conjecture in terms of the multilinear version of the Fourier representation \eqref{lambdafrequency}; the multilinear version of \eqref{lambdaspace} will look somewhat differently, involving a power of the determinant as well as a signum function in case of even powers of the determinant.
\begin{conjecture}
For $n\geq 4$ set
\begin{equation}
\Lambda_n(f_1,\dots,f_n) :=
\textup{p.v.} \int_{(\mathbb{R}^{n-1})^n} \bigg(\prod_{j=1}^n\widehat{f_j}(\xi_j )\bigg)
\det\begin{pmatrix}1 & \!\cdots\! & 1\\ \xi_1 & \!\cdots\! & \xi_n \end{pmatrix}^{-1} d\xi_1\cdots d\xi_n.
\end{equation}
If
\begin{equation}
n-1<p_j<\infty \text{ for } j=1,\ldots,n \ \text{ and } \ \sum_{j=1}^n \frac{1}{p_j} = 1,
\end{equation}
then
\begin{equation}
|\Lambda_n(f_1,\dots,f_n)| \lesssim_{p_1,\ldots,p_n} \prod_{j=1}^n \|f_j\|_{\textup{L}^{p_j}(\mathbb{R}^{n-1})}.
\end{equation}
\end{conjecture}
Our approach to Theorem~\ref{triltheorem} does not suffice to prove this conjecture without further ideas.

\section{Convergence of the principal value}
\label{convergence}

\subsection{Existence of the principal value in \eqref{lambdaspace}}
In this subsection we prove that \eqref{lambdaspace} is well-defined. Along the lines of the same arguments one can prove that \eqref{lambdafrequency}
is also well-defined; in fact our reasoning follows the ideas by Valdimarsson from \cite{Val}, who has given an explicit absolutely integrable expression for \eqref{lambdafrequency}.

Our goal is to prove that the limit
\begin{equation}\label{limitepsilon}
\lim_{\epsilon\to 0^+} \int_{|\det| > \epsilon}\limits \frac{f(x)g(y)h(z)}{\det\begin{pmatrix}1 & 1 & 1\\ x & y & z \end{pmatrix}} \delta(x+y+z)\, dxdydz
\end{equation}
exists for any three Schwartz functions $f,g,h$. Here we have abbreviated
$$ \det=\det\begin{pmatrix}1 & 1 & 1\\ x & y & z \end{pmatrix} $$
in the specification of the domain of integration. The above integral may be viewed as the Lebesgue integral over a subset of $\mathbb{R}^4$ by introducing orthonormal coordinates on a four dimensional subspace of $\mathbb{R}^6$ determined by
$$ x+y+z=0 $$
and rewriting the integral in these coordinates. Then we also need to divide the obtained expression by $3$, accounting for two factors of $\sqrt{3}$, the length of the vector $(1,1,1)$, taken into account twice since the Dirac delta in \eqref{limitepsilon} is two dimensional.

If we denote $\varphi(x,y,z) = f(x)g(y)h(z)$, we only need to observe that $\varphi$ is now a Schwartz function on $\mathbb{R}^6$.
It is clear that the integral in \eqref{limitepsilon} is absolutely integrable for a fixed $\epsilon>0$ and it can be written as an iterated integral
\begin{equation}\label{iterated}
\int_{\mathbb{R}^3} \Bigg[ \,\int_{|\det| > \epsilon}\limits \frac{\varphi(x,y,z)}{\det\begin{pmatrix}1 & 1 & 1\\ x & y & z \end{pmatrix}}
\delta(x_1+y_1+z_1)\, dx_1dy_1dz_1 \Bigg]\delta(x_2+y_2+z_2)\, dx_2dy_2dz_2.
\end{equation}
We may restrict attention to both $(x_1,y_1,z_1)$ and $(x_2,y_2,z_2)$ perpendicular to $(1,1,1)$ and nonzero, which is the case on the domain of integration, and parameterize
$$ (x_1,y_1,z_1)=\lambda e + \mu e^\perp, $$
where $e=(e_1,e_2,e_3)$ is the unit vector in the direction of $(x_2,y_2,z_2)$ and $e^\perp=(e_1^\perp,e_2^\perp,e_3^\perp)$ is the unit vector in the direction of the cross product $(x_2,y_2,z_2)\times (1,1,1)$. Then the value of the determinant in \eqref{iterated} becomes simply $\mu \varrho \sqrt{3}$ with $\varrho =(x_2^2+y_2^2+z_2^2)^{1/2}$.
Define
\begin{align*}
x^\ast &= (\lambda e_1 - \mu e_1^\perp,x_2),\\
y^\ast &= (\lambda e_2 -\mu e_2^\perp,y_2),\\
z^\ast &= (\lambda e_3 -\mu e_3^\perp,z_2),
\end{align*}
so that we have
$$ {\det\begin{pmatrix}1 & 1 & 1\\ x^\ast & y^\ast & z^\ast \end{pmatrix}} = -{\det\begin{pmatrix}1 & 1 & 1\\ x & y & z \end{pmatrix}}. $$
We may rewrite \eqref{iterated} using symmetrization relative to the reflection across the origin of the innermost integration variable $\mu$:
\begin{equation}\label{symmetrized}
\frac{1}{3} \int_{\mathbb{R}^3} \int_{\mathbb{R}} \bigg[ \int_{|\mu|>\epsilon/(\varrho \sqrt{3})} \big(\varphi(x,y,z)-\varphi(x^\ast,y^\ast,z^\ast)\big) \frac{d\mu}{2\mu} \bigg]
d\lambda \frac{\delta(x_2+y_2+z_2)}{\varrho }\, dx_2dy_2dz_2.
\end{equation}
For fixed $x_2,y_2,z_2,\lambda$ the difference $\varphi(x,y,z)-\varphi(x^\ast,y^\ast,z^\ast)$ is smooth in $\mu$, rapidly decaying as $|\mu|\to\infty$, and vanishes at $\mu=0$. Hence
\begin{align}
& \lim_{\epsilon\to 0^+} \int_{|\mu|>\epsilon/(\varrho \sqrt{3})} \big(\varphi(x,y,z)-\varphi(x^\ast,y^\ast,z^\ast)\big) \frac{d\mu}{2\mu} \nonumber \\
& = \int_{\mathbb{R}} \big(\varphi(x,y,z)-\varphi(x^\ast,y^\ast,z^\ast)\big) \frac{d\mu}{2\mu} \label{pointwise}
\end{align}
in the sense of the Lebesgue integral on the right hand side, and in particular the limit on the left hand side exists. In order to see that the whole expression \eqref{symmetrized} converges as $\epsilon\to 0^+$, we use the Lebesgue dominated convergence theorem. Note that for any $(x,y,z)$ satisfying
$$ |(x,y,z) - (x^\ast,y^\ast,z^\ast)| = 2|\mu| \leq 1 $$
the Mean Value Theorem gives an intermediate point $(x^{\ast\ast},y^{\ast\ast},z^{\ast\ast})$ such that
$$ \big|\varphi(x,y,z)-\varphi(x^\ast,y^\ast,z^\ast)\big| \leq 2|\mu| |\nabla\varphi(x^{\ast\ast},y^{\ast\ast},z^{\ast\ast})| \lesssim_{\varphi} 2|\mu|, $$
which combined with the Schwartz decay of $\varphi$ and $\nabla\varphi$ yields
$$ \left| \frac{\varphi(x,y,z)-\varphi(x^\ast,y^\ast,z^\ast)}{2\mu} \right| \lesssim_{\varphi} (1+\mu^2)^{-1} (1+\lambda^2)^{-1} (1+\varrho ^2)^{-1}. $$
Finally, in the outer integral we pass to the polar coordinates on $x_2+y_2+z_2=0$ and note the finiteness of
$$ \int_{0}^{2\pi}\! \int_{0}^{\infty} \Big(\int_\mathbb{R} \frac{d\mu}{1+\mu^2}\Big) \Big(\int_\mathbb{R} \frac{d\lambda}{1+\lambda^2}\Big) \frac{1}{(1+\varrho^2)\varrho } \,\varrho  d\varrho  d\theta. $$
Hence convergence of \eqref{symmetrized} as $\epsilon \to 0^+$ follows from \eqref{pointwise} by the dominated convergence theorem. For later purpose we note that we have also shown that the inner integral in \eqref{iterated} converges pointwise almost everywhere as $\epsilon\to 0^+$, i.e.\@ if $x_2+y_2+z_2 = 0$, then there exists
\begin{equation}\label{limitepsilon2}
\lim_{\epsilon\to 0^+} \int_{\{(x_1,y_1,z_1)\in\mathbb{R}^3 : |\det| > \epsilon\}}\limits \frac{f(x)g(y)h(z)}{\det\begin{pmatrix}1 & 1 & 1\\ x & y & z \end{pmatrix}} \delta(x_1+y_1+z_1)\, dx_1dy_1dz_1,
\end{equation}
and the limiting function integrates on this plane to the value of \eqref{limitepsilon}.

\subsection{Equality between \eqref{lambdaspace} and \eqref{lambdafrequency}}
Let us split $\mathbb{R}^6$ orthogonally (in skew coordinates) as $\mathbb{R}^4\times \mathbb{R}^2$, where $\mathbb{R}^4$ coordinatizes $x+y+z=0$ and $\mathbb{R}^2$ provides coordinates for its orthogonal complement.  More precisely, let $e_1$, $e_2$, $e_3$ be a positively oriented orthonormal basis of $\mathbb{R}^3$, with $e_3 = \frac{1}{\sqrt{3}} (1,1,1)$. Let $(\tilde{x},\tilde{y},\tilde{z}) \in (\mathbb{R}^2)^3$ be the components of $(x,y,z)$ with respect to the frame $e_1,e_2,e_3$. In other words, $(\tilde{x},\tilde{y},\tilde{z})$ are defined by
\begin{equation} \label{eq:tildevar}
\left( \begin{array}{ccc} 
\tilde{x} & \tilde{y} & \tilde{z}
\end{array} \right) R
=  \left( \begin{array}{ccc} 
{x} & {y} & {z}
\end{array} \right),
\end{equation}
where $R$ is the rotation matrix whose $i$-th row is $e_i$, for $i = 1,2,3$. Then
$$
x+y+z = \sqrt{3}\tilde{z},
$$
and
$$
\det \left( \begin{array}{ccc} 
1 & 1 & 1 \\
{x} & {y} & {z}
\end{array} \right)
= 
\det 
\left[ \left( \begin{array}{ccc} 
0 & 0 & \sqrt{3} \\
\tilde{x} & \tilde{y} & \tilde{z}
\end{array} \right) R \right] = \sqrt{3} \det \left( \begin{array}{cc} 
\tilde{x} & \tilde{y} 
\end{array} \right) ,
$$
so (1.1) can be rewritten
\begin{equation} \label{eq:det2paper}
\frac{1}{3\sqrt{3}} \int_{(\mathbb{R}^2)^3} f(x) g(y) h(z) \delta(\tilde{z}) \det \left( 
\begin{array}{cc}
\tilde{x} & \tilde{y}
\end{array}
\right)^{-1} d\tilde{x} d\tilde{y} d\tilde{z},
\end{equation}
where $(x,y,z)$ are defined by (\ref{eq:tildevar}). Now let $T$ be the distribution on $\mathbb{R}^4$, defined by
$$
T\phi = \textrm{p.v.} \int_{\mathbb{R}^4} \phi(\tilde{x},\tilde{y}) \det(\tilde{x} \, \tilde{y})^{-1} d\tilde{x} d\tilde{y} = \lim_{\epsilon \to 0^+} \int_{|\det(x \, y)| > \epsilon} \phi(x,y) \det(x \, y)^{-1} dx dy.
$$
We claim that the Fourier transform of $T$ is itself. Then the Fourier transform of the distribution $\delta(\tilde{z}) \otimes T$ is $1 \otimes T$. Also, write $(\xi,\eta,\zeta)$ for the dual variables of $(x,y,z)$, and $(\tilde{\xi},\tilde{\eta},\tilde{\zeta})$ for the dual variables to $(\tilde{x},\tilde{y},\tilde{z})$, so that
\begin{equation} \label{eq:tildevardual}
\left( \begin{array}{ccc} 
\tilde{\xi} & \tilde{\eta} & \tilde{\zeta}
\end{array} \right) R
=  \left( \begin{array}{ccc} 
{\xi} & {\eta} & {\zeta}
\end{array} \right).
\end{equation}
Then the Fourier transform of $f(x)g(y)h(z)$, as a function of $(\tilde{x},\tilde{y},\tilde{z})$, is just $\widehat{f}(\xi) \widehat{g}(\eta) \widehat{h}(\zeta)$. Hence (\ref{eq:det2paper}) can be rewritten as
\begin{equation} \label{eq:det3paper}
\frac{1}{3\sqrt{3}} \int_{(\mathbb{R}^2)^3} \widehat{f}(\xi) \widehat{g}(\eta) \widehat{h}(\zeta) \det \left( 
\begin{array}{cc}
\tilde{\xi} & \tilde{\eta}
\end{array}
\right)^{-1} d\tilde{\xi} d\tilde{\eta} d\tilde{\zeta}.
\end{equation}
But $\det  \left( 
\begin{array}{cc}
\tilde{\xi} & \tilde{\eta}
\end{array}
\right)$ can be rewritten as
$$
\frac{1}{\sqrt{3}} \det
\left( \begin{array}{ccc}
0 & 0 & \sqrt{3} \\ 
\tilde{\xi} & \tilde{\eta} & \tilde{\zeta}
\end{array} \right) 
= \frac{1}{\sqrt{3}} \det\left[ \left( \begin{array}{ccc} 
1 & 1 & 1 \\
{\xi} & {\eta} & {\zeta}
\end{array} \right) R^{-1} \right]
= \frac{1}{\sqrt{3}} \det \left( \begin{array}{ccc} 
1 & 1 & 1 \\
{\xi} & {\eta} & {\zeta}
\end{array} \right).
$$
Hence (\ref{eq:det3paper}) is equal to (\ref{lambdafrequency}); it remains to prove our claim above, which we proceed as follows.

As in the previous subsection one can show
\begin{equation}\label{chainbeginning}
T(\phi) = \int_{\mathbb{R}^2} \Bigg[\lim_{\epsilon\to 0^+} \int_{\{y\in\mathbb{R}^2:|\det(x\,y)|>\epsilon\}}\limits \!\!\!\phi(x,y) \det(x\ y)^{-1} dy \Bigg] dx.
\end{equation}
In the inner integral we perform a linear, orthonormal, symmetric, and involutory change of variables:
$$ z=(z_1,z_2)=Ry:=\left(\frac{x\cdot y}{|x|},\frac{\det(x\ y)}{|x|}\right). $$
Setting $\tilde{\phi}(x,z)=\phi(x,y)$ we obtain
$$ T(\phi) = \int_{\mathbb{R}^2} \bigg[ \lim_{\epsilon\to 0^+} \int_{\{z\in\mathbb{R}^2:|z_2|>\epsilon/|x|\}}\limits \!\!\!\tilde{\phi}(x,z) |x|^{-1}z_2^{-1} dz \bigg] dx. $$
Up to the normalization factor $|x|^{-1}$, the inner integral is a pairing of the Schwartz function $\tilde{\phi}(x,.)$ with the distribution $\mathbbm{1}\otimes \mathop{\textup{p.v.}}1/t$, which is the Fourier transform of $\delta\otimes i\pi\,\textup{sgn}$. Hence, denoting the partial Fourier transform in the $j$-th variable by $\mathcal{F}_j$, we get
$$ T(\phi) = i\pi \int_{\mathbb{R}^2} \left[ \int_{\mathbb{R}^2} \mathcal{F}_2 {\tilde{\phi}}(x,\zeta) \delta(\zeta_1)
|x|^{-1} \mathop{\textup{sgn}}\zeta_2 \,d\zeta \right] dx, $$
so substituting
$$ \zeta=(\zeta_1,\zeta_2)=R\eta=\left(\frac{x\cdot\eta}{|x|},\frac{\det(x\ \eta)}{|x|}\right) $$
yields
\begin{align}
T(\phi) & = i\pi \int_{\mathbb{R}^2} \left[ \int_{\mathbb{R}^2} \mathcal{F}_2 {\phi}(x,\eta) \delta\Big(\frac{x\cdot\eta}{|x|}\Big)
|x|^{-1} \mathop{\textup{sgn}}\Big(\frac{\det(x\ \eta)}{|x|}\Big) d\eta \right] dx \nonumber \\
& = i\pi \int_{\mathbb{R}^4}  \mathcal{F}_2 {\phi}(x,\eta) \delta({x\cdot\eta})
\mathop{\textup{sgn}}(\det(x\ \eta)) \,d\eta dx. \label{chainend}
\end{align}
Observing that $\mathcal{F}_1^2$ is the reflection about the origin in the first argument and using invariance of the whole integral under simultaneous reflection in both variables we obtain
\begin{align*}
T(\phi) & = i\pi \int_{\mathbb{R}^4}  \mathcal{F}_1^2\mathcal{F}_2 {\phi}(x,-\eta) \delta({x\cdot \eta})
\mathop{\textup{sgn}}(\det(x \,\eta)) d\eta dx \\
& = i\pi \int_{\mathbb{R}^4}  \mathcal{F}_1 \widehat{\phi}(x,\eta) \delta({\eta \cdot x})
\mathop{\textup{sgn}}(\det(\eta \,x)) dx d\eta.
\end{align*}
In the last identity we changed variable $\eta \mapsto -\eta$, and used that $\det(x \,\, -\eta) = \det(\eta \,\, x)$.
Applying the symmetric argument of the above chain of identities leading from \eqref{chainbeginning} to \eqref{chainend}, the last integral expression simplifies to $T(\widehat{\phi})$. This proves $T=\widehat{T}$ and thus also establishes the equality of \eqref{lambdaspace} and \eqref{lambdafrequency}.

\section{Proof of the main theorem}
\label{mainproof}

\subsection{Proof of the a priori estimate \eqref{eqformbound1}}
Fix a point $(x_2,y_2,z_2)$ perpendicular to $(1,1,1)$ and such that none of its components vanish. The change of variables
$$ \tilde{x}_1 = x_1/x_2,\ \ \tilde{y}_1 = y_1/y_2,\ \ \tilde{z}_1 = z_1/z_2 $$
and the new functions
$$ {\setlength{\arraycolsep}{2pt}\begin{array}{rll}
\tilde{f}(\tilde{x}_1) & := (\mathop{\textup{sgn}}x_2) \,f(x_1,x_2) & = (\mathop{\textup{sgn}}x_2) \,f(\tilde{x}_1 x_2,x_2), \\[1.5mm]
\tilde{g}(\tilde{y}_1) & := (\mathop{\textup{sgn}}y_2) \,g(y_1,y_2) & = (\mathop{\textup{sgn}}y_2) \,g(\tilde{y}_1 y_2,y_2), \\[1.5mm]
\tilde{h}(\tilde{z}_1) & := (\mathop{\textup{sgn}}z_2) \,h(z_1,z_2) & = (\mathop{\textup{sgn}}z_2) \,h(\tilde{z}_1 z_2,z_2)
\end{array}} $$
allow to rewrite the integral in \eqref{limitepsilon2} as
\begin{align*}
& \int_{|\det| > \epsilon}\limits \tilde{f}(\tilde{x}_1)\tilde{g}(\tilde{y}_1)\tilde{h}(\tilde{z}_1) \,\det\begin{pmatrix}1 & 1 & 1\\ x_1 & y_1 & z_1 \\ x_2 & y_2 & z_2 \end{pmatrix}^{-1} \delta(x_1+y_1+z_1)\, x_2y_2z_2\, d\tilde{x}_1d\tilde{y}_1d\tilde{z}_1 \nonumber \\
& = \int_{|\det| > \epsilon}\limits \tilde{f}(\tilde{x}_1)\tilde{g}(\tilde{y}_1)\tilde{h}(\tilde{z}_1)
\,\det\begin{pmatrix}x_2^{-1} & y_2^{-1} & z_2^{-1}\\ \tilde{x}_1 & \tilde{y}_1 & \tilde{z}_1 \\ 1 & 1 & 1 \end{pmatrix}^{-1}
\delta(x_2\tilde{x}_1+y_2\tilde{y}_1+z_2\tilde{z}_1)\, d\tilde{x}_1d\tilde{y}_1d\tilde{z}_1. \label{changed}
\end{align*}
The limit of this expression as $\epsilon\to 0^+$ becomes \eqref{bht} with $(a,b,c)=(x_2,y_2,z_2)$. By a simple scaling argument using $\delta(\alpha x) =|\alpha|^{-1}\delta(x)$ we may replace $(a,b,c)$ with a unit vector in the same direction without changing the value of \eqref{bht}.

Now we intend to apply the main result of Grafakos and Li \cite{GL}. In order to relate \eqref{bht} to the more common definition of the bilinear Hilbert transform we introduce coordinates $(s,t)$ on the plane $ax+by+cz=0$ such that
$$ (x,y,z) = s(1,1,1) + t(u,v,w), $$
where $(u,v,w)$ is the unit vector in the direction of $(1,1,1)\times (a,b,c)$. Since
$$ \begin{pmatrix} 1/\sqrt{3} & a & u \\ 1/\sqrt{3} & b & v \\ 1/\sqrt{3} & c & w \end{pmatrix} $$
is an orthogonal matrix with determinant $1$, we can calculate
\begin{align*}
& \det\begin{pmatrix}1 & 1 & 1\\ a^{-1} & b^{-1} & c^{-1} \\s+ut & s+vt & s+wt \end{pmatrix}
= t \,\det\begin{pmatrix}1 & 1 & 1\\ a^{-1} & b^{-1} & c^{-1} \\u & v & w \end{pmatrix} \\
& = \frac{t}{\sqrt{3}} \det\left[ \begin{pmatrix}1 & 1 & 1\\ a^{-1} & b^{-1} & c^{-1} \\u & v & w \end{pmatrix}
\begin{pmatrix} 1 & a & u \\ 1 & b & v \\ 1 & c & w \end{pmatrix} \right]
= \frac{t}{\sqrt{3}} \det\begin{pmatrix}3 & 0 & 0\\ \ast & 3 & \ast\\ 0 & 0 & 1 \end{pmatrix} = 3\sqrt{3}\,t,
\end{align*}
so that \eqref{bht} turns into
\begin{equation}\label{dualbht}
\frac{1}{3} \,\textup{p.v.} \int_{\mathbb{R}^2} \tilde{f}(s+ut) \tilde{g}(s+vt) \tilde{h}(s+wt) \,\frac{dt}{t} \,ds.
\end{equation}
By a trivial change of variables
$$ \tilde{s}=s+wt,\ \  \tilde{t}=(u-w)t $$
and denoting $\kappa=(v-w)/(u-w)$ we see that \eqref{dualbht} becomes
$$ \frac{1}{3} \int_\mathbb{R} \Big(\textup{p.v.} \int_{\mathbb{R}} \tilde{f}(\tilde{s}+\tilde{t}) \tilde{g}(\tilde{s}+\kappa\tilde{t}) \,\frac{d\tilde{t}}{\tilde{t}}\Big) \,\tilde{h}(\tilde{s}) d\tilde{s}, $$
which in turn can be identified as the trilinear form corresponding to the bilinear operator estimated in \cite{GL}. For the purpose of the discussion in Section~\ref{remarks} we also note that indeed
$$ \kappa= \frac{v-w}{u-w} = -\frac{a}{b} = -\frac{x_2}{y_2}. $$

Applying the main result from the paper \cite{GL} we can bound the absolute value of \eqref{limitepsilon2} by a constant multiple of
\begin{align*}
& \|\tilde{f}\|_{\textup{L}^p(\mathbb{R})}\|\tilde{g}\|_{\textup{L}^q(\mathbb{R})}\|\tilde{h}\|_{\textup{L}^r(\mathbb{R})} \\
 = & |x_2|^{-1/p} \|f(\cdot,x_2)\|_{\textup{L}^p(\mathbb{R})}\, |y_2|^{-1/q} \|g(\cdot,y_2)\|_{\textup{L}^q(\mathbb{R})}\, |z_2|^{-1/r} \|h(\cdot,z_2)\|_{\textup{L}^r(\mathbb{R})}
\end{align*}
for any choice of exponents $p,q,r$ as in Theorem~\ref{triltheorem}. Define
$$ F(x_2) := \|f(\cdot,x_2)\|_{\textup{L}^p(\mathbb{R})},\ \ G(y_2) := \|g(\cdot,y_2)\|_{\textup{L}^q(\mathbb{R})},\ \ H(z_2) := \|h(\cdot,z_2)\|_{\textup{L}^r(\mathbb{R})}. $$
Using what we showed in Section~\ref{convergence} we can now estimate \eqref{limitepsilon} by a constant depending on $p,q,r$ times
$$ \int_{\mathbb{R}^3} \frac{F(x_2)G(y_2) H(z_2)}{|x_2|^{1/p} |y_2|^{1/q} |z_2|^{1/r}} \,\delta(x_2+y_2+z_2)\, dx_2 dy_2 dz_2. $$
Let us introduce polar coordinates in the plane $x_2+y_2+z_2=0$. Thus we write
$$ (x_2,y_2,z_2)= \varrho (e_1(\theta),e_2(\theta),e_3(\theta)) $$
with $\varrho=(x_2^2+y_2^2+z_2^2)^{1/2}\in (0,\infty)$ and $\theta\in [0,2\pi)$. The last integral is
$$ \frac{1}{\sqrt{3}} \int_{0}^{2\pi}\! \int_0^{\infty} \frac{F(\varrho e_1(\theta))G(\varrho e_2(\theta))H(\varrho e_3(\theta))}{|e_1(\theta)|^{1/p} |e_2(\theta)|^{1/q} |e_3(\theta)|^{1/r}} \,d\varrho d\theta $$
and then using H\"older's inequality in the inner integration it is further reduced to
$$ \int_{0}^{2\pi } \frac{\|F\|_{\textup{L}^p(\mathbb{R})} \|G\|_{\textup{L}^q(\mathbb{R})} \|H\|_{\textup{L}^r(\mathbb{R})}}{|e_1(\theta)|^{2/p} |e_2(\theta)|^{2/q} |e_3(\theta)|^{2/r}} \,d\theta. $$
It remains to observe that
$$ \int_{0}^{2\pi} |e_1(\theta)|^{-2/p} |e_2(\theta)|^{-2/q} |e_3(\theta)|^{-2/r} d\theta <\infty, $$
since the singularities $e_i(\theta)=0$ occur at six different angles $\theta$ and are each integrable due to the conditions $p,q,r>2$. Thus inequality \eqref{eqformbound1} follows from $\|F\|_{\textup{L}^p(\mathbb{R})}=\|f\|_{\textup{L}^p(\mathbb{R}^2)}$, $\|G\|_{\textup{L}^q(\mathbb{R})}=\|g\|_{\textup{L}^q(\mathbb{R}^2)}$, and $\|H\|_{\textup{L}^p(\mathbb{R})}=\|h\|_{\textup{L}^p(\mathbb{R}^2)}$.

\subsection{Sharpness of the range}
Take $1\leq p,q,r\leq \infty$. We will assume $1/p+1/q+1/r=1$, for otherwise a simple scaling argument proves that no bound of the form \eqref{eqformbound1} can hold. Suppose that we do not have $2<p,q,r<\infty$. By symmetry and interpolation it then suffices to assume $p\leq 2$ and that $q$ and $r$ are finite.

Fix perpendicular unit vectors $(x_1^0,y_1^0,z_1^0)$ and $(x_2^0,y_2^0,z_2^0)$, both perpendicular to $(1,1,1)$ and with
$$ \det\begin{pmatrix}1/\sqrt{3} & 1/\sqrt{3} & 1/\sqrt{3}\\ x^0 & y^0 & z^0 \end{pmatrix} = 1. $$
Note that $|x^0|=|y^0|=|z^0|=\sqrt{2/3}$, since the above matrix is orthogonal.
Also fix $\delta=1/100$ and consider some $\varrho\geq 1$. If $x$ is in $\textup{B}_{2\delta}(x^0)$, the ball of radius $2\delta$ about $x^0$, and if $y\in \textup{B}_{2\varrho\delta}(\varrho y^0)$ and $z\in \textup{B}_{2\varrho\delta}(-x^0-\varrho y^0)$ are such that $x+y+z=0$, then
$$ \det\begin{pmatrix}1 & 1 & 1\\ x & y & z \end{pmatrix} = 3 \det(x\ y) = 3\varrho \det\big(x\ y/\varrho\big), $$
while similarly
$$ \sqrt{3} = \det\begin{pmatrix}1 & 1 & 1\\ x^0 & y^0 & z^0 \end{pmatrix} = 3 \det(x^0\ y^0). $$
Since $x$ is a $2\delta$-perturbation of $x^0$ and $y/\varrho$ is a $2\delta$-perturbation of $y^0$, we have
\begin{equation}\label{eqdetest}
0<\det\begin{pmatrix}1 & 1 & 1\\ x & y & z \end{pmatrix}\leq 2\varrho.
\end{equation}

Take a $\textup{C}^\infty$ function $\phi$ with values in $[0,1]$, which is supported on $\textup{B}_2(0,0)$ and constantly equals $1$ on $\textup{B}_1(0,0)$. For any $t>0$ and $c\in\mathbb{R}^2$ denote
$$ \phi_{t,c}(x) := \phi\Big(\frac{x-c}{t}\Big). $$
Simply let $f:=\phi_{\delta, x^0}$. For a positive integer $N\geq 10$ set $g:=\sum_{n=10}^N g_n$, where
$$ g_n:=2^{-2n/q}n^{-2/q} \phi_{2^n\delta, 2^ny^0}, $$
and $h:=\sum_{n=10}^N h_n$, where
$$ h_n:=2^{-2n/r}n^{-2/r} \phi_{2^n\delta, -x^0-2^ny^0}. $$
We clearly have $\|f\|_{\textup{L}^p(\mathbb{R}^2)}<\infty$, as well as (by the disjointness of the supports thanks to $|y^0|>8\delta$)
$$ \|g\|_{\textup{L}^q(\mathbb{R}^2)}^q = \sum_{n=10}^N  \|g_n\|_{\textup{L}^q(\mathbb{R}^2)}^q\lesssim \sum_{n=10}^N 2^{-2n}n^{-2} 2^{2n} \lesssim 1, $$
and similarly $\|h\|_{\textup{L}^r(\mathbb{R}^2)}^r \lesssim 1$ independently of $N$.
On the other hand, a simple support argument together with \eqref{eqdetest} for $\varrho=2^n$ shows
\begin{align*}
\Lambda(f,g,h)
& = \sum_{n=10}^N \int_{(\mathbb{R}^2)^3} \frac{f(x)g_n(y)h_n(z)}{\det\begin{pmatrix}1 & 1 & 1\\ x & y & z \end{pmatrix}} \,\delta(x+y+z)\, dxdydz \\
& \geq \sum_{n=10}^N 2^{-n-1} \!\!\int_{\substack{|x-x^0|<\delta\\ |y-2^n y^0|<2^n\delta\\ |z+x^0+2^n y^0|<2^n\delta}}\limits\!\!
2^{-2n/q}n^{-2/q} 2^{-2n/r}n^{-2/r} \,\delta(x+y+z)\, dxdydz \\
& \gtrsim \sum_{n=10}^N 2^{-n} 2^{-2n/q}n^{-2/q} 2^{-2n/r}n^{-2/r} 2^{2n} = \sum_{n=10}^N 2^{(2/p-1)n} n^{2/p-2}.
\end{align*}
Because of $p\leq 2$ the last sum grows unboundedly with $N$, proving that there is no a priori bound of the form \eqref{eqformbound1} for the considered case of exponents.

Alternatively, a slightly more explicit counterexample can be constructed from appropriate smooth truncations of the functions
\begin{align*}
f(x_1,x_2) & = x_1^{-2/p} (\log x_1)^{-2/p}, \\
g(y_1,y_2) & = y_1^{-2/q} (\log y_1)^{-2/q}, \\
h(z_1,z_2) & = (-z_1)^{-2/r} \big(\log(-z_1)\big)^{-2/r}.
\end{align*}

\section{Closing remarks}
\label{remarks}

\subsection{Calder\'{o}n's commutator}
The first commutator of Calder\'{o}n is defined by
\begin{align*}
\mathcal{C}(F,G)(x) & := \textup{p.v.} \int_{\mathbb{R}} F(y) \frac{G(x)-G(y)}{(x-y)^2} dy \\
& \,= \textup{p.v.} \int_{0}^{1}\! \int_{\mathbb{R}} F(y) \,G'((1-\kappa) x + \kappa y)\, \frac{dy}{x-y} \,d\kappa
\end{align*}
and substituting $t=y-x$ we can write it as
$$ \mathcal{C}(F,G)(x) = - \,\textup{p.v.} \int_{0}^{1}\! \int_{\mathbb{R}} F(x+t) G'(x+\kappa t) \,\frac{dt}{t} \,d\kappa. $$
Suppose that we want to deduce the estimate
\begin{equation}\label{caldcommest}
\|\mathcal{C}(F,G)\|_{\textup{L}^{r'}(\mathbb{R})} \lesssim_{p,q,r} \|F\|_{\textup{L}^{p}(\mathbb{R})} \|G'\|_{\textup{L}^{q}(\mathbb{R})}
\end{equation}
for exponents $p,q,r$ in the region \eqref{eqexpreg} solely from Theorem~\ref{triltheorem}. If we formally substitute
\begin{align*}
f(x_1,x_2) & := -2^{1/p} \mathbbm{1}_{(-1,0)}(x_2) F(x_1/x_2), \\
g(y_1,y_2) & := 2^{1/q} \mathbbm{1}_{(0,1)}(y_2) G'(y_1/y_2), \\
h(z_1,z_2) & := -2^{1/r} \mathbbm{1}_{(-1,0)}(z_2) H(z_1/z_2)
\end{align*}
into $\Lambda(f,g,h)$, then the computation presented in Section~\ref{mainproof} gives
\begin{align*}
& \frac{2}{3\sqrt{3}} \int_{\mathbb{R}^2} \Big( \textup{p.v.} \int_{\mathbb{R}^2} F(\tilde{s}+\tilde{t}) G'(\tilde{s}+\kappa\tilde{t}) H(\tilde{s}) \,\frac{d\tilde{t}}{\tilde{t}} \,d\tilde{s} \Big) \\
& \qquad\qquad \mathbbm{1}_{(-1,0)}(-\kappa y_2) \mathbbm{1}_{(0,1)}(y_2) \mathbbm{1}_{(-1,0)}((\kappa-1)y_2) \,|y_2| \,dy_2 d\kappa,
\end{align*}
which is up to an unimportant constant equal to
$$ \int_{0}^{1} \textup{p.v.} \int_{\mathbb{R}^2} F(\tilde{s}+\tilde{t}) G'(\tilde{s}+\kappa\tilde{t}) H(\tilde{s}) \,\frac{d\tilde{t}}{\tilde{t}}\, d\tilde{s} \,d\kappa
= - \int_{\mathbb{R}} \mathcal{C}(F,G)(x) H(x) dx. $$
Since we also have
$$ \|f\|_{\textup{L}^p(\mathbb{R}^2)} = \|F\|_{\textup{L}^p(\mathbb{R})},\ \ \|g\|_{\textup{L}^q(\mathbb{R}^2)} = \|G'\|_{\textup{L}^q(\mathbb{R})},\ \ \|h\|_{\textup{L}^r(\mathbb{R}^2)} = \|H\|_{\textup{L}^r(\mathbb{R})}, $$
the established bound \eqref{eqformbound1} becomes just dualized version of the desired estimate \eqref{caldcommest}. Of course, these functions $f,g,h$ are not in the Schwartz space (and thus not in the domain of $\Lambda$), so indeed one also has to apply standard approximation arguments.

\subsection{Mixed norm estimates}
Observe that the dilation symmetry group of $\Lambda$, which is $\textup{GL}_2(\mathbb{R})$, particularly contains non-isotropic dilations,
$$ \bigg\{ \begin{pmatrix}\alpha & 0\\ 0 & \beta\end{pmatrix} \,:\, \alpha,\beta>0 \bigg\}. $$
That way $\Lambda$ can also be viewed as a bi-parameter object and mixed $\textup{L}^p$ norms,
$$ \|f\|_{\textup{L}^{p_2}(\textup{L}^{p_1})} := \big\| \|f(x_1,x_2)\|_{\textup{L}^{p_1}_{x_1}} \big\|_{\textup{L}^{p_2}_{x_2}}
= \left[ \Big( \int_{\mathbb{R}} |f(x_1,x_2)|^{p_1} dx_1 \Big)^{p_2/p_1} dx_2 \right]^{1/p_2}, $$
come up naturally in that setting. Exactly the same proof gives
$$ |\Lambda(f,g,h)| \lesssim_{p_1,p_2,q_1,q_2,r_1,r_2} \|f\|_{\textup{L}^{p_2}(\textup{L}^{p_1})} \,\|g\|_{\textup{L}^{q_2}(\textup{L}^{q_1})} \,\|h\|_{\textup{L}^{r_2}(\textup{L}^{r_1})} $$
for six finite exponents satisfying
$$ \frac{1}{p_i}+\frac{1}{q_i}+\frac{1}{r_i}=1;\ \ i=1,2,\quad \frac{1}{p_1}+\frac{1}{p_2}<1,\quad \frac{1}{q_1}+\frac{1}{q_2}<1,\quad \frac{1}{r_1}+\frac{1}{r_2}<1, $$
as long as we have at our disposal an $\textup{L}^{p_1}\times\textup{L}^{q_1}\times\textup{L}^{r_1}$ bound for \eqref{bht}. Unfortunately, the largest possible range of estimates for the bilinear Hilbert transform is not known.

Let us conclude by recalling that the bi-parameter bilinear Hilbert transform from \cite{MPTT} does not satisfy any $\textup{L}^p$ estimates, which makes the positive result of Theorem~\ref{triltheorem} seem less expected.

\section*{Acknowledgments}
The authors would like to thank the American Institute of Mathematics for hosting the workshop \emph{Carleson theorems and multilinear operators} in May 2015. The collaboration was initiated and the majority of the work was performed as a part of the workshop program.
P. G. was also partially supported by the NSF grant DMS-1361697 and an Alfred P. Sloan Research Fellowship.
D. H. was also partially supported by the Simons Foundation Grant Number 315380.
V. K. was also partially supported by the Croatian Science Foundation under the project 3526.
B. S. was also partially supported by the NSF grant DMS-1401671.
C. T. was also partially supported by the Hausdorff Center for Mathematics.
P.-L. Y. was also partially supported by a direct grant for research from the Chinese University of Hong Kong (4053120).

\end{document}